\documentclass[12pt]{article}

\title{Edge Growth in Graph Cubes\thanks{This research was done at the Graph Coloring workshop at the Technion, Haifa Israel.}}

\author{
   Matt DeVos\thanks{Supported in part by an NSERC Discovery Grant (Canada) and a Sloan Fellowship.}
 \and
   St\'ephan Thomass\'e\thanks{Universit\'e Montpellier 2 - CNRS, LIRMM 161 rue Ada, 34392 Montpellier
    Cedex, France
    \tt{thomasse@lirmm.fr}}}
   
\date{}

\usepackage{fullpage}
\usepackage{graphicx}
\usepackage{amsthm}
\usepackage{epsfig}
\usepackage{amsfonts}
\usepackage{amsmath}
\usepackage{enumitem}
\usepackage{mathrsfs}

\theoremstyle{plain}
\newtheorem{theorem}{Theorem}[section]

\newtheorem{conjecture}[theorem]{Conjecture}

\theoremstyle{definition}

\begin{document}

\maketitle

\begin{abstract}
We show that for every connected graph $G$ of diameter $\ge 3$, the graph $G^3$ has average degree $\ge \frac{7}{4} \delta(G)$.  We also provide an example showing that this bound is best possible.  This resolves a question of Hegarty~\cite{PH}.
\end{abstract}

\section{Introduction}

Throughout the paper, we only consider simple graphs. Let $G$ be a graph. We denote 
by $v(G)$, $e(G)$ its number of vertices, edges respectively, and let $\delta(G)$ denote the minimum degree of $G$. The 
{\em $k^{th}$-power} of $G$, denoted by $G^k$, has vertex set $V(G)$ and edges the pair of 
vertices at distance at most $k$ in $G$. If $G$ is connected, the \emph{diameter} of $G$ 
is the maximum distance between a pair of vertices of $G$, or, equivalently, 
the smallest integer $k$ so that $G^k$ is a clique.

Consider a generating set $A$ of a finite (multiplicative) group and suppose that $1 \in A$ and $g \in A \Rightarrow g^{-1} \in A$.  Numerous important questions in Number Theory and Group Theory concern the increase in size from $|A|$ to $|A^k|$.  Such problems can be phrased naturally in terms of Cayley graphs.  If $G$ is the (simple) Cayley graph generated by $A$, then $G^k$ is generated by $A^k$ and the sizes of the sets $A$ and $A^k$ are given by the degrees of these (regular) graphs.  Thus the growth of the set $A^k$ can be studied in terms of the number of additional edges in the graph $G^k$.  For instance, the following result is an easy corollary of a famous theorem of Cauchy and Davenport.

\begin{theorem}[Cauchy-Davenport]
If $G$ is a connected Cayley graph on a group of prime order with diameter $<k$ then
$e(G^k) \ge k e(G)$.
\end{theorem}

Inspired by this connection, Hegarty considered the more general problem of how many extra edges are formed when we move from a graph $G$ to the $k^{th}$ power of $G$.  Although little can be said for graphs in general, the problem is interesting for connected regular graphs with a diameter constraint.  Perhaps surprisingly, even for this class of graphs, there does not exist a positive constant $c$ so that $e(G^2) \ge (1+c) e(G)$.  In contrast to this, the following holds for the third power:

\begin{theorem}[Hegarty]
There exists a positive constant $c$ so that every connected regular graph of diameter $\ge 3$ satisfies $e(G^3) \ge (1+c) e(G)$.
\end{theorem}

Hegarty proved this for $c=0.087$ and this was subsequently improved by Pokrovskiy~\cite{AP} who showed 
that the same result holds with $c= \frac{1}{6}$ (Pokrovskiy also established some results for higher powers of $G$).  These authors both raised the question of the best possible value of $c$.  We settle this problem in the following theorem.

\begin{theorem} 
\label{main}
If $G$ is a connected graph with diameter $\ge 3$, then $e(G^3) \ge \frac{7}{8} \delta(G) v(G)$.
\end{theorem}

In particular, when $G$ is regular, this shows that $c$ can be chosen to be $\frac{3}{4}$.
To see that this is best possible, we construct a family of regular graphs defined as follows.  
The graph $G_k$ is obtained from the disjoint union of the graphs $H_1,H_2,\ldots,H_5$ by adding all possible edges between vertices in $H_i$ and $H_{i+1}$ for $1 \le i \le 4$, where the graphs $H_1$ and $H_5$ are copies of $K_{2k+1}$, the graphs $H_2$ and $H_4$ are copies of $K_{2k}$ minus a perfect matching, and $H_3$ is a single vertex.  
It follows 
that $G_k$ is $4k$-regular with $8k+3$ vertices so $e(G_k) = \frac{1}{2}(8k+3)(4k) = 16k^2 + 6k$.  Its cube $G_k^3$ 
has $4k+1$ vertices of degree $8k+2$ and $4k+2$ vertices of degree $6k+1$ so it satisfies 
$e(G_k^3) = \frac{1}{2}(4k+1)(8k+2) + \frac{1}{2}(4k+2)(6k+1) = 28k^2 + 16k + 2$.  The family of 
graphs $\{G_k\}_{k \in {\mathbb N}}$ hence shows that the constant $\frac{7}{8}$ in Theorem~\ref{main} is best possible.

\begin{figure}[h]
\centering
\includegraphics[height=3cm]{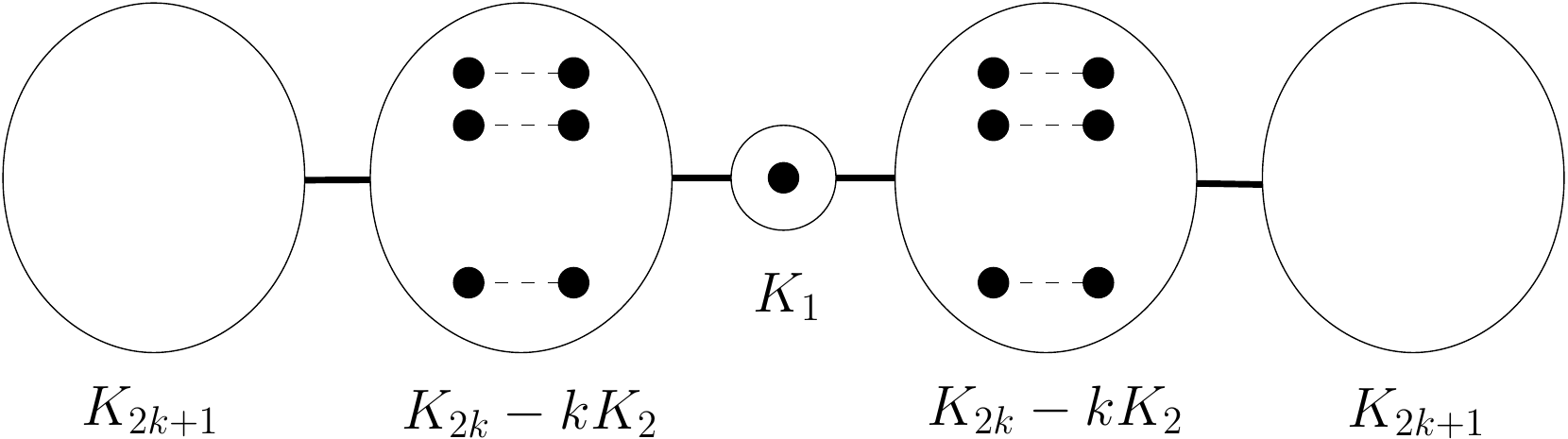}
\caption{The graph $G_k$}
\end{figure}

There are a number of interesting related problems for directed graphs.  Here we highlight a rather basic conjecture, which, if true, would resolve a 
special case of the Caccetta-H\"aggkvist conjecture.

\begin{conjecture}
If $D$ is an orientation of a simple graph and every vertex of $D$ has indegree and outdegree equal to $d$ then $e(D^2) \ge 2e(D)$.
\end{conjecture}

\section{Proof}

For a set of vertices $X$ we let $N(X)$ denote the closed neighbourhood of $X$, i.e. $N(X)$ is the union of $X$ and the set of vertices with a neighbour in $X$.  For a nonnegative integer $k$ we let $N^k(X)$ denote the set of vertices at distance $\leq k$ from a point in $X$.  For a vertex $v$ we simplify this notation by $N(v) = N( \{v\})$ and $N^k(v) = N^k(\{v\})$.  Note that the degree of a $v$ in $G^3$ satisfies ${\mathit deg}_{G^3}(v) = |N^3(v)| - 1$.

\bigskip

\noindent{\it Proof of Theorem \ref{main}:} 
Let $G$ be a connected graph with minimum degree $\delta$ and diameter $\ge 3$.  
We say that a path is \emph{geodesic} if it is a shortest path between its endpoints.  A vertex $v$ is \emph{doubling} if ${\mathit deg}_{G^3}(v) \ge 2 \delta$. We let $Z$ be the set of doubling vertices in $G$. We now prove a sequence of claims. 

\bigskip

\noindent{(1) If $v$ is an internal vertex in a geodesic path of length 3, then $v$ is doubling.

\smallskip

To see this, suppose that our geodesic path has vertex sequence $u,v,v',u'$.  Now $N(u) \cap N(u') = \emptyset$ and $N(u) \cup N(u') \subseteq N^3(v)$ so $v$ is doubling.

\bigskip

Now let $X_1, X_2, \ldots, X_m$ be the vertex sets of the components of $G - Z$.  

\bigskip

\noindent{(2) If $v$ and $v'$ both belong to the same $X_i$, for some $1 \le i \le m$, then $N^2(v) = N^2(v')$.

\smallskip

Since $G[X_i]$ is connected, it suffices to prove that $N^2(v) \subseteq N^2(v')$ when $v,v'$ are adjacent.  In this case, suppose that $u \in N^2(v)$.  Then there is a path of length 3 from $v'$ to $u$ which has $v$ 
as an internal vertex.  By (1) this path cannot be geodesic, so there must be a path of length at most
$2$ from $v'$ to $u$, i.e. $u \in N^2(v')$.  

\bigskip

Next, define a relation $\sim$ on $\{X_1, \ldots, X_m \}$ by the rule that $X_i \sim X_j$ if 
$N(X_i) \cap N(X_j) \neq \emptyset$.

\bigskip

\noindent{(3) If $X_i \sim X_j$, $v \in X_i$ and $v' \in X_j$, then 
$N^2(v) = N^2(v')$.

\smallskip

In light of (2), it suffices to prove this in the case that $N(v) \cap N(v') \neq \emptyset$.  
To see this, suppose (for a contradiction) that $u \in N^2(v) \setminus N^2(v')$.  Then we have
$N(u) \cap N(v') = \emptyset$ and $N(u) \cup N(v') \subseteq N^3(v)$ so $v$ is doubling, which is contradictory.  

\bigskip

\noindent{(4)} $\sim$ is an equivalence relation.

\smallskip

To check that $\sim$ is transitive, suppose that $X_i \sim X_j \sim X_k$ and choose 
$v \in X_i$ and $v' \in X_k$.  It follows from (3) that $N^2(v) = N^2(v')$ but then $v$ and $v'$ 
have a common neighbour, hence $N(X_i) \cap N(X_k) \neq \emptyset$.

\bigskip

Let $\{Y_1,Y_2,\ldots,Y_{\ell}\}$ be the set of unions of equivalence classes of $\sim$.
\bigskip

\noindent{(5)} The subgraph of $G^2$ induced by $N(Y_i)$ is a clique for every $1 \le i \le \ell$.

\smallskip

Let $v,v' \in N(Y_i)$.  If one of $v,v'$ is in $Y_i$ then it follows from (3) that $v,v'$ are adjacent in $G^2$.  In the remaining case, choose $u \in Y_i$ adjacent to $v$.  Since $v' \in N^2(u)$ there is a path of length $\le 3$ from $v$ to $v'$ which has $u$ as an internal vertex.  It now follows from (1) that $v$ and $v'$ are distance $\le 2$ in $G$, so they are adjacent in $G^2$.

\bigskip

Let $y_i = |Y_i|$ for every $1 \le i \le \ell$.

\bigskip

\noindent{(6)} ${\mathit deg}_{G^3}(v) \ge \delta + y_i$ for every $v \in Y_i$.

\smallskip

Claim (5) shows that $N(Y_i)$ induces a clique in $G^2$.  Since $G$ has diameter $\ge 3$ the graph $G^2$ is not a clique.  Hence there must exist a vertex 
$u \in N^2(Y_i) \setminus N(Y_i)$.  Now $N(u) \cap Y_i = \emptyset$ and $N(u) \cup Y_i \subseteq N^3(v)$ which gives us ${\mathit deg}_{G^3}(v) \ge \delta + y_i$ as desired.

\bigskip

Set $y = y_1 + y_2 + \ldots + y_{\ell}$ and set $z = |Z|$.  

\bigskip

\noindent{(7) $z \ge \delta \ell - y$

\smallskip


First note that $\delta \le |N(Y_i)| = |Y_i| + |N(Y_i) \cap Z|$ so $|N(Y_i) \cap Z| \ge \delta - y_i$.  Next, observe that $N(Y_i) \cap N(Y_j) = \emptyset$ whenever $i \neq j$.  This gives us $z = |Z| \ge \sum_{i = 1}^{\ell} |N(Y_i) \cap Z| \ge \sum_{i=1}^{\ell} (\delta - y_i) = \delta \ell - y$ as desired.

\bigskip

We now have the tools to complete the proof.  Combining the fact that every vertex in $Z$ has degree at least $2 \delta$ in $G^3$ with (6), 
gives us the following inequality (here we use Cauchy-Schwarz and (7) in getting to the third line)
\begin{align*}
\sum_{v \in V(G)} {\mathit deg}_{G^3}(v) - \tfrac{7}{4} \delta v(G) 
        &\ge 2 \delta z + \sum_{i = 1}^{\ell} y_i (\delta + y_i) - \tfrac{7}{4} \delta (z+y) \\
        &= \tfrac{1}{4}\delta z - \tfrac{3}{4} \delta y + \sum_{i=1}^{\ell} y_i^2 \\
        &\ge \tfrac{1}{4} \delta (\delta \ell - y) - \tfrac{3}{4} \delta y + \frac{y^2}{\ell} \\
        &= \left( \frac{\delta \sqrt{\ell}}{2} - \frac{y}{\sqrt{\ell}} \right)^2 \\
        & \ge 0.
\end{align*}
This shows that $G^3$ has average degree $\ge \frac{7}{4} \delta$, thus completing the proof.
\quad\quad$\Box$

\end{document}